\input amstex
\magnification=1200
\documentstyle{amsppt}
\NoRunningHeads
\NoBlackBoxes
\topmatter
\title Interactive games and representation theory
\endtitle
\author\bf D.V.Juriev
\endauthor
\date math.FA/9803020\enddate
\affil\rm
``Thalassa Aitheria'' Research Center for Mathematical Physics
and Informatics,\linebreak
ul.Miklukho-Maklaya, 20-180, Moscow 117437, Russia.\linebreak
E-mail: denis\@juriev.msk.ru.
\endaffil
\abstract This short note is a conceptual prologue to the series of articles
devoted to the analysis of interrelations between the representation theory
(especially, its inverse problems) and the control and games theory.
\endabstract
\endtopmatter
\document
This note was appeared as an attempt to formulate rigorously and
ma\-the\-ma\-ti\-cal\-ly the ideas of [1] and simultaneously to treat
them in a more (mathematically) general context of control and games
theories (see e.g.[2,3]).

\definition{Definition 1} {\it An interactive system\/} (with $n$
{\it interactive controls\/}) is a control system with $n$ independent controls
coupled with unknown or incompletely known feedbacks (the feedbacks as well
as their couplings with controls are of a so comp\-li\-cated nature that their
can not be described completely). {\it An interactive game\/} is a game with
interactive controls of each player.
\enddefinition

Below we shall consider only deterministic and differential interactive
systems. For symplicity we suppose that $n=2$. In this case the general
interactive system may be written in the form:
$$\dot\varphi=\Phi(\varphi,u_1,u_2),\tag1$$
where $\varphi$ characterizes the state of the system and $u_i$ are
the interactive controls:
$$u_i(t)=u_i(u_i^\circ(t),\left.[\varphi(\tau)]\right|_{\tau\le t}),$$
i.e. the independent controls $u_i^\circ(t)$ coupled with the feedbacks on
$\left.[\varphi(\tau)]\right|_{\tau\le t}$. One may suppose that the
feedbacks are integrodifferential on $t$.

\proclaim{Theorem} Each interactive system (1) may be transformed to the
form (2) below (which is not, however, unique):
$$\dot\varphi=\tilde\Phi(\varphi,\xi),\tag2$$
where the magnitude $\xi$ (with infinite degrees of freedom as a rule)
obeys the equation
$$\dot\xi=\Xi(\xi,\varphi,\tilde u_1,\tilde u_2),\tag3$$
where $\tilde u_i$ are the interactive controls of the form $\tilde
u_i(t)=\tilde u_i(u_i^\circ(t); \varphi(t),\xi(t))$ (i.e. the feedbacks
are on $\xi(t)$ as well as on $\varphi(t)$ and are differential on $t$).
\endproclaim

\remark{Remark 1} One may exclude $\varphi(t)$ from the feedbacks in
the interactive controls $\tilde u_i(t)$.
\endremark

\definition{Definition 2} The magnitude $\xi$ with its dynamical equations
(3) and its cont\-ri\-bution into the interactive controls $\tilde u_i$ will
be called {\it the intention field}.
\enddefinition

The theorem is the formal mathematical counterpart for the main dynamical
hypothesis of [1].

\remark{Remark 2} The theorem holds true for the interactive games.
\endremark

\definition{Definition 3} {\it The main inverse problem of representation
theory for the in\-ter\-ac\-tive system (1)\/} (or for the interactive
game) is
\roster
\item"(1)" to write the system (1) in the form (2);
\item"(2)" to determine the geometrical and algebraical structure of
the intention field;
\item"(3)" to find the algebraic structure, which ``governs'' the dynamics
(3), e.g. solving the dynamical inverse problem of representation theory
[4] (see also [5]) for the system (3).
\endroster
\enddefinition

\remark{Remark 3} The solution of the main inverse problem of representation
theory for the interactive system may use {\sl a posteriori} data on the
system.
\endremark

The main inverse problem of representation theory for the interactive systems
will be one of the topics in the series of articles [6] devoted to the novel
interrelations between the representation and the control theories.

\remark{Remark 4} Visualizing the intention field(s) somehow one is able to
enlarge the initial controlled system and to strengthen its controllability
considering the controls $u^\circ(t)$ depending on ``positional'' analysis
of $\xi$ as well as of $\varphi$.
\endremark

\Refs
\roster
\item" [1]" Juriev D., On the dynamics of physical interactive information
systems. E-print: mp\_arc/97-158.
\item" [2]" Pontryagin L.S. et al., Mathematical theory of optimal processes.
Moscow, 1969 [in Russian]; Differential geometric control theory.
Eds.R.W.Brockett, R.S.Millman, H.J.Suss\-mann, Birkhauser, 1983.
\item" [3]" Isaacs R., Differential games. A mathematical theory with
applications to war\-fare and pur\-suit, control and optimization. Wiley, New
York, 1965; Owen G., Games theory. Saun\-ders, Philadelphia, 1968;
Vorob'ev N.N., Uspekhi Matem. Nauk 15(2) (1970) 80-140 [in Russian].
\item" [4]" Juriev D., Dynamical inverse problem of representation theory and
noncommutative geometry [in Russian]. Fundam.Prikl.Matem. 4(1) (1998)
[e-version: funct-an/9507001 (1995)].
\item" [5]" Juriev D., An excursus into the inverse problem of representation
theory [in Russian]. Report RCMPI-95/04 (1995) [e-version: mp\_arc/96-477
(1996)].
\item" [6]" Juriev D., Topics in hidden algebraic structures and
infinite-dimensional dynamical sym\-met\-ries of controlled systems, in
preparation.
\endroster
\endRefs
\enddocument